\documentclass[11pt]{article}
\usepackage{epic,latexsym,amssymb}
\usepackage{tikz}

\usepackage{tikz,epsf}

\usepackage{amsfonts}
\usepackage{amscd}
\usepackage{amsmath}
\usepackage{graphicx}
\usepackage{color}
\usepackage{caption,subcaption}

\newtheorem{thm}{Theorem}

\newtheorem{ob}[thm]{Observation}
\newtheorem{prop}[thm]{Proposition}
\newtheorem{lem}[thm]{Lemma}

\newtheorem{cor}[thm]{Corollary}

\newcommand{\ecc}{{\rm ecc}}
\newcommand{\diam}{{\rm diam}}
\newcommand{\rad}{{\rm rad}}

\newcommand{\proof}{\noindent\textbf{Proof. }}

\newcommand{\qed}{$\Box$}

\newcommand{\1}{\vspace{0.1cm}}

\newcommand{\vertex}{\node[vertex]}
\tikzstyle{vertex}=[circle, draw, inner sep=0pt, minimum size=6pt]

\newcommand{\QEDmark}{\mbox{\textsc{qed}}}
\newcommand{\proofStarter}[1]{\textsc{#1} }

\def\vertex(#1){\put(#1){\circle*{2}}}
\def\vertexo(#1){\put(#1){\circle{2}}}
\def\vert(#1){\put(#1){\circle*{1.5}}}
\def\verto(#1){\put(#1){\circle{1.5}}}
\def\lab(#1)#2{\put(#1){\makebox(0,0)[c]{#2}}}
\begin{document}

\title{Lower Bounds on the Distance Domination Number \\ of a Graph}
\author{$^{1,3}$Randy Davila, $^2$Caleb Fast, $^1$Michael A. Henning\thanks{Research supported in part by the South African National Research Foundation and the University of Johannesburg}\, and  $^2$Franklin Kenter \\
\\
$^1$Department of Mathematics\\
University of Johannesburg \\
Auckland Park 2006, South Africa \\
\small {\tt Email: mahenning@uj.ac.za} \\
\\
$^2$Computational and Applied Mathematics \\
Rice University \\
Houston, TX 77005, USA \\
\small {\tt Email: calebfast@gmail.com, fhk2@rice.edu}\\
\\
$^3$Department of Mathematics \\
Texas State University-San Marcos \\
San Marcos, TX 78666, USA \\
\small {\tt Email: randyrdavila@gmail.com}
}

\date{}
\maketitle

\begin{abstract}
For an integer $k \ge 1$, a (distance) $k$-dominating set of a connected graph $G$ is a set $S$ of vertices of $G$ such that  every vertex of $V(G) \setminus S$ is at distance at most~$k$ from some vertex of $S$. The $k$-domination number, $\gamma_k(G)$, of $G$ is the minimum cardinality of a $k$-dominating set  of $G$. In this paper, we establish lower bounds on the $k$-domination number of a graph in terms of its diameter, radius and girth. We prove that for connected graphs $G$ and $H$, $\gamma_k(G \times H) \ge \gamma_k(G) + \gamma_k(H) -1$, where $G \times H$ denotes the direct product of $G$ and $H$.
\end{abstract}


\section{Introduction}

Distance in graphs is a fundamental concept in graph theory. Let $G$ be a connected graph. The \emph{distance} between two vertices $u$ and $v$ in $G$, denoted $d_G(u,v)$,
is the length (i.e., the number of edges) of a shortest $(u,v)$-path in $G$. The \emph{eccentricity} $\ecc_G(v)$ of $v$ in $G$ is the distance between $v$ and a vertex farthest from $v$ in $G$. The minimum eccentricity among all vertices of $G$ is the \emph{radius} of $G$, denoted by $\rad(G)$, while the maximum eccentricity among all vertices of $G$ is the \emph{diameter} of $G$, denoted by $\diam(G)$. Thus, the diameter of $G$ is the maximum distance among all pairs of vertices of $G$.
A vertex $v$ with $\ecc_G(v) = \diam(G)$ is called a \emph{peripheral vertex} of $G$.
A \emph{diametral path} in $G$ is a shortest path in $G$ whose length is equal to the diameter of the graph. Thus, a diametral path is a path of length $\diam(G)$ joining two peripheral vertices of $G$.
If $S$ is a set of vertices in $G$, then the \emph{distance}, $d_G(v,S)$, from a vertex $v$ to the set $S$ is the minimum distance from $v$ to a vertex of $S$; that is, $d_G(v,S) = \min \{d_G(u,v) \mid u \in S\}$. In particular, if $v \in S$, then $d(v,S) = 0$.

Domination in graphs is also very well studied in graph theory. A \emph{dominating set} in a graph $G$ is a set $S$ of vertices of $G$ such that every vertex in $V(G) \setminus S$ is adjacent to at least one vertex in $S$. The \emph{domination number} of $G$, denoted by $\gamma(G)$, is the minimum cardinality of a dominating set of $G$. The literature on the subject of domination parameters in graphs up to the year 1997 has been surveyed and detailed in the two books~\cite{hhs1, hhs2}.

In this paper, we continue the study of distance domination in graphs which combines the concepts of both distance and domination in graphs. Let $k \ge 1$ be an integer and let $G$ be a graph.  In 1975, Meir and Moon~\cite{MeMo75} introduced the concept of a distance $k$-dominating set (called a ``$k$-covering" in~\cite{MeMo75}) in a graph. A set $S$ is a $k$-\emph{dominating set} of $G$ if every vertex is within distance~$k$ from some vertex of $S$; that is, for every vertex $v$ of $G$, we have $d(v,S) \le k$.  The $k$-\emph{domination number} of $G$, denoted  $\gamma_k(G)$, is the minimum cardinality of a $k$-dominating set of $G$. When $k = 1$, the $1$-domination number of $G$ is precisely the domination number of $G$; that is, $\gamma_1(G) = \gamma(G)$.
The literature on the subject of distance domination in graphs up to the year 1997 can be found in the book chapter~\cite{He98}. Distance domination is now widely studied, see, for example,~\cite{CyLeRa,Fr88,HaMeVo07,HeLi15+,HeOeSw,MeVo05,MeMo75,Sl76,TiXu09,ToVo91}.

\medskip
\noindent\textbf{Definitions and Notation.}
For notation and graph theory terminology, we in general follow~\cite{MHAYbookTD}.
Specifically, let $G$ be a graph with vertex set
$V(G)$ of order~$n(G) = |V(G)|$ and edge set $E(G)$ of size~$m(G) = |E(G)|$.
We assume throughout the paper that all graphs considered are \emph{simple} graphs, i.e., finite graphs with no directed edges and no loops.
A \emph{non-trivial graph} is a graph on at least two vertices.
A \emph{neighbor} of a vertex $v$ in $G$ is a vertex adjacent to $v$. The \emph{open neighborhood} of $v$, denoted $N_G(v)$, is the set of all neighbors of $v$ in $G$, while the \emph{closed  neighborhood} of $v$ is the set $N_G[v] = N_G(v) \cup \{v\}$.
The \emph{closed $k$-neighborhood}, denoted $N_k[v]$, of $v$ is defined in~\cite{Fr88} as the set of all vertices within distance~$k$ from $v$ in $G$; that is, $N_k[v] = \{u \mid d(u,v) \le k\}$. When $k = 1$, the set $N_k[v] = N[v]$.

The \emph{degree} of a vertex $v$ in $G$, denoted $d_G(v)$, is the number of neighbors, $|N_G(v)|$, of $v$ in $G$. The minimum and maximum degree among all the vertices of $G$ are denoted by $\delta(G)$ and $\Delta(G)$, respectively.
The subgraph induced by a set $S$ of vertices of $G$ is denoted by $G[S]$.
The \emph{girth} of $G$, denoted $g(G)$, is the length of a shortest cycle in $G$.
For sets of vertices  $X$ and $Y$ of $G$, the set $X$ $k$-\emph{dominates} the set $Y$ if every vertex of $Y$ is within distance~$k$ from some vertex of $X$. In particular, if $X$ $k$-dominates the set $V(G)$, then $X$ is a $k$-dominating set of $G$.

If the graph $G$ is clear from context, we simply write $V$, $E$, $d(v)$, $\ecc(v)$, $N(v)$ and $N[v]$ rather than $V(G)$, $E(G)$, $d_G(v)$, $\ecc_G(v)$, $N_G(v)$ and $N_G[v]$, respectively. We use the standard notation $[n] = \{1,2,\ldots,n\}$.

\medskip
\noindent\textbf{Known Results.}
The $k$-domination number of $G$ is in the class of $NP$-hard graph invariants to compute~\cite{hhs2}. Because of the computational complexity of computing $\gamma_k(G)$, graph theorists have sought upper and lower bounds on $\gamma_k(G)$ in terms of simple graph parameters like order, size, and degree.

In 1975, Meir and Moon~\cite{MeMo75} established an upper bound for the $k$-domination number of a tree in terms of its order. They proved that for $k \ge 1$, if $T$ is a tree of order~$n \ge k+1$, then $\gamma_k(T) \le n/(k + 1)$. As a consequence of this result and Observation~\ref{ob:spanning}, if $G$ is a connected graph of order~$n \ge k + 1$, then
$\gamma_k(G) \le \frac{n}{k + 1}$. A short proof of the Meir-Moon upper bound can also be found in~\cite{HeOeSw} (see, also, Proposition~24 and Corollary~12.5 in the book chapter~\cite{He98}).
A complete characterization of the graphs $G$ achieving equality in this upper bound was obtained by Topp and Volkmann~\cite{ToVo91}.
Tian and Xu~\cite{TiXu09} improved the Meir-Moon upper bound and showed that for $k \ge 1$, if $G$ is a connected graph of order~$n \ge k + 1$ with maximum degree~$\Delta$, then $\gamma_k(G) \le \frac{1}{k}(n - \Delta + k - 1)$.
The Tian-Xu bound was further improved by Henning and Lichiardopol~\cite{HeLi15+} who showed that for $k \ge 2$, if $G$ is a connected graph with minimum degree~$\delta \ge 2$ and maximum degree~$\Delta$ and of order~$n \ge \Delta + k - 1$, then $\gamma_k(G) \le \frac{n + \delta - \Delta}{\delta + k - 1}$.

We recall the following well-known lower bound on the domination number of a graph in terms of its diameter.

\begin{thm}{\rm (\cite{hhs2})}
If $G$ is a connected graph with diameter~$d$, then $\gamma(G) \ge \frac{d+1}{3}$.
\label{thm:diam}
\end{thm}

The following two results were originally conjectured by the conjecture making program Graffiti.pc (see~\cite{DeLa}).

\begin{thm}{\rm (\cite{DePeWa10})} \label{thm:radius}
If $G$ is a connected graph with radius~$r$, then $\gamma(G) \ge \frac{2r}{3}$.
\end{thm}

\begin{thm}{\rm (\cite{DePeWa10})}
\label{thm:girth}
If $G$ is a connected graph with girth $g \ge 3$, then $\gamma(G) \ge \frac{g}{3}$.
\end{thm}

\noindent\textbf{Our Results.}  In this paper, we establish lower bounds for the $k$-domination number of a graph in terms of its diameter (Theorem~\ref{thm:Diameter}), radius (Corollary~\ref{cor:Radius}), and girth (Theorem~\ref{thm:Girth}). These results generalize the results of Theorem~\ref{thm:diam}, Theorem~\ref{thm:radius}, and Theorem~\ref{thm:girth}. A key tool in order to prove our results is the important lemma (Lemma~\ref{ob:spanning}) that every connected graph has a spanning tree with equal $k$-domination number. We also prove a key property (Lemma~\ref{Lem:CycleOutsider}) of shortest cycles in a graph that enables us to establish our girth result for the $k$-domination number of a graph. We show that our bounds are all sharp and examples are provided following the proofs.

\section{Preliminary Observations and Lemmas}

Since every $k$-dominating set of a spanning subgraph of a graph $G$ is a $k$-dominating set of $G$, we have the following observation.

\begin{ob}
For $k \ge 1$, if $H$ is a spanning subgraph of a graph $G$, then $\gamma_k(G) \le \gamma_k(H)$.
\label{ob:spanning}
\end{ob}

We shall also need the following lemma.

\begin{lem}
\label{ob:spanning}
For $k \ge 1$, every connected graph $G$ has a spanning tree $T$ such that $\gamma_k(T) = \gamma_k(G)$.
\label{lem:spanning}
\end{lem}
\proof Let $S = \{v_1, \ldots, v_\ell\}$ be a minimum $k$-dominating set of $G$. Thus, $|S| = \ell = \gamma_k(G)$. We now partition the vertex set $V(G)$ into $\ell$ sets $V_1,\ldots,V_\ell$ as follows. Initially, we let $V_i = \{v_i\}$ for all $i \in [\ell]$. We then consider sequentially the vertices not in $S$. For each vertex $v \in V(G) \setminus S$, we select a vertex $v_i \in S$ at minimum distance from~$v$ in $G$ and add the vertex $v$ to the set $V_i$. We note that if $v \in V(G) \setminus S$ and $v \in V_i$ for some $i \in [\ell]$, then $d_G(v,v_i) = d_G(v,S)$, although the vertex $v_i$ is not necessarily the  unique vertex of $S$ at minimum distance from~$v$ in $G$. Further, since $S$ is a $k$-dominating set of $G$, we note that $d_G(v,v_i) \le k$. For each $i \in [\ell]$, let $T_i$ be a spanning tree of $G[V_i]$ that is distance preserving from the vertex~$v_i$; that is, $V(T_i) = V_i$ and for every vertex $v \in V(T_i)$, we have $d_{T_i}(v,v_i) = d_G(v,v_i)$. We now let $T$ be the spanning tree of $G$ obtained from the disjoint union of the $\ell$ trees $T_1,\ldots,T_\ell$ by adding $\ell - 1$ edges of $G$. We remark that these added $\ell - 1$ edges exist as $G$ is connected. We now consider an arbitrary vertex, $v$ say, of $G$. The vertex $v \in V_i$ for some $i \in [\ell]$. Thus, $d_T(v,v_i) \le d_{T_i}(v,v_i) = d_G(v,v_i) = d_G(v,S) \le k$. Therefore, the set $S$ is a $k$-dominating set of $T$, and so $\gamma_k(T) \le |S| = \gamma_k(G)$. However, by Observation~\ref{ob:spanning}, $\gamma_k(G) \le \gamma_k(T)$. Consequently, $\gamma_k(T) = \gamma_k(G)$.~\qed

\begin{lem}
\label{Lem:CycleOutsider}
Let $G$ be a connected graph that contains a cycle, and let $C$ be a shortest cycle in $G$. If $v$ is a vertex of $G$ outside $C$ that $k$-dominates at least $2k$ vertices of $C$, then there exist two vertices $u,w \in V(C)$ that are both $k$-dominated by $v$ and such that a shortest $(u,v)$-path does not contain $w$ and a shortest $(v,w)$-path does not contain $u$.
\end{lem}
\proof Since $v$ is not on $C$, it has a distance of at least~$1$ to every vertex of $C$.  Let $u$ be a vertex of $C$ at minimum distance from~$v$ in $G$. Let $Q$ be the set of vertices on $C$ that are $k$-dominated by $v$ in $G$. Thus, $Q \subseteq V(C)$ and, by assumption, $|Q| \ge 2k$. Among all vertices in $Q$, let $w \in Q$ be chosen to have maximum distance from~$u$ on the cycle $C$. Since there are $2k-1$ vertices within distance~$k-1$ from $u$ on $C$, the vertex $w$ has distance at least~$k$ from $u$ on the cycle $C$. Let $P_u$ be a shortest $(u,v)$-path and let $P_w$ be a shortest $(v,w)$-path in $G$. If $w \in V(P_u)$, then $d_G(v,w) < d_G(v,u)$, contradicting our choice of the vertex~$u$. Therefore, $w \notin V(P_u)$. Suppose that $u \in V(P_w)$. Since $C$ is a shortest cycle in $G$, the distance between $u$ and $w$ on $C$ is the same as the distance between $u$ and $w$ in $G$. Thus, $d_G(u,w) = d_C(u,w)$, implying that $d_G(v,w) = d_G(v,u) + d_G(u,w) \ge 1 + d_G(u
 ,w) = 1 + d_C(u,w) \ge 1 + k$, a contradiction.  Therefore, $u \notin V(P_w)$.~\qed

\section{Lower Bounds}

In this section we provide various lower bounds on the $k$-domination number for general graphs. We first prove a generalization of Theorem~\ref{thm:diam} by establishing a lower bound on the $k$-domination number of a graph in terms of its diameter. We remark that when $k = 1$, Theorem~\ref{thm:Diameter} is precisely Theorem~\ref{thm:diam}.

\begin{thm}
\label{thm:Diameter}
For $k \ge 1$, if $G$ is a connected graph with diameter $d$, then
\[
\gamma_{k}(G) \ge \frac{d+1}{2k+1}.
\]
\end{thm}
\proof Let $P \colon u_0u_1 \ldots u_d$ be a diametral path in $G$, joining two peripheral vertices $u = u_0$ and $v = u_d$ of $G$. Thus, $P$ has length~$\diam(G) = d$. We show that every vertex of $G$ $k$-dominates at most~$2k+1$ vertices of $P$. Suppose, to the contrary, that there exists a vertex $q \in V(G)$ that $k$-dominates at least $2k+2$ vertices of $P$. (Possibly, $q \in V(P)$.) Let $Q$ be the set of vertices on the path $P$ that are $k$-dominated by the vertex~$q$ in $G$. By supposition, $|Q| \ge 2k+2$. Let $i$ and $j$ be the smallest and largest integers, respectively, such that $u_i \in Q$ and $u_j \in Q$. We note that $Q \subseteq \{u_i,u_{i+1},\ldots,u_j\}$. Thus, $2k+2 \le |Q| \le j - i + 1$. Since $P$ is a shortest $(u,v)$-path in $G$, we therefore note that $d_G(u_i,u_j) = d_P(u_i,u_j) = j - i \ge 2k+1$.
Let $P_i$ be a shortest $(u,q)$-path in $G$ and let $P_j$ be a shortest $(q,v)$-path in $G$. Since the vertex $q$ $k$-dominates both $u_i$ and $u_j$ in $G$, both paths $P_u$ and $P_v$ have length at most~$k$. Therefore, the $(u_i,u_j)$-path obtained by following the path $P_i$ from $u_i$ to $q$, and then proceeding along the path $P_j$ from $q$ to $u_j$, has length at most~$2k$, implying that $d_G(u_i,u_j) \le 2k$, a contradiction. Therefore, every vertex of $G$ $k$-dominates at most~$2k+1$ vertices of $P$.

Let $S$ be a minimum $k$-dominating set of $G$. Thus, $|S| = \gamma_k(G)$. Each vertex of $S$ $k$-dominates at most~$2k+1$ vertices of $P$, and so $S$ $k$-dominates at most~$|S|(2k+1)$ vertices of $P$. However, since $S$ is a  $k$-dominating set of $G$, every vertex of $P$ is $k$-dominated the set $S$, and so $S$ $k$-dominates $|V(P)| = d + 1$ vertices of $P$. Therefore, $|S|(2k+1) \ge d + 1$, or, equivalently, $\gamma_k(G) = (d+1)/(2k+1)$.~\qed

\medskip
That the lower bound of Theorem~\ref{thm:Diameter} is tight may be seen by taking $G$ to be a path, $v_1v_2 \ldots v_n$, of order~$n = \ell (2k+1)$ for some $\ell \ge 1$. Let $d = \diam(G)$, and so $d = n-1 =  \ell (2k+1) - 1$. By Theorem~\ref{thm:Diameter}, $\gamma_k(G) \ge (d + 1)/(2k+1) = \ell$. The set
\[
S = \bigcup_{i=0}^{\ell - 1} \{v_{k+1 + i(2k+1)}\}
\]
is a $k$-dominating set of $G$, and so $\gamma_k(G) \le |S| = \ell$. Consequently, $\gamma_k(G) = \ell = (d + 1)/(2k+1)$. We state this formally as follows.

\begin{prop}
If $G = P_n$ where $n \equiv 0 \bmod (2k+1)$, then $\gamma_k(G) = \frac{\diam(G) + 1}{2k+1}$.
\label{prop:path}
\end{prop}

More generally, by applying Theorem~\ref{thm:Diameter}, the $k$-domination number of a cycle $C_n$ or path $P_n$ on $n \ge 3$ vertices is easy to compute.

\begin{prop}
For $k \ge 1$ and $n \ge 3$,
$\gamma_k(P_n) = \gamma_k(C_n) = \lceil \frac{n}{2k+1} \rceil$.
\label{prop:cycle}
\end{prop}

By replacing each vertex $v_i$, for $2 \le i \le n -1$, on the path $v_1v_2 \ldots v_n$ with a clique (clique $V_i$ corresponds to vertex $v_i$) of size at least~$\delta \ge 1$, and adding all edges between $v_1$ and vertices in $V_2$, adding all edges between $v_n$ and vertices in $V_{n-1}$, and adding all edges between vertices in $V_i$ and $V_{i+1}$ for $2 \le i \le n -2$, we obtain a graph with minimum degree~$\delta$ achieving the lower bound of Theorem~\ref{thm:Diameter}.

As a consequence of Theorem~\ref{thm:Diameter}, we have the following lower bound on the $k$-domination number of a graph in terms of its radius. We remark that when $k = 1$, Corollary~\ref{cor:Radius} is precisely Theorem~\ref{thm:radius}. Therefore, Corollary~\ref{cor:Radius} is a generalization of Theorem~\ref{thm:radius}.

\begin{cor}
\label{cor:Radius}
For $k \ge 1$, if $G$ is a connected graph with radius~$r$, then
\[
\gamma_{k}(G) \ge \frac{2r}{2k+1}.
\]
\end{cor}
\proof
By Lemma~\ref{lem:spanning}, the graph $G$ has a spanning tree $T$ such that $\gamma_k(T) = \gamma_k(G)$. Since adding edges to a graph cannot increase its radius, $\rad(G) \le \rad(T)$. Since $T$ is a tree, we note that $\diam(T) \ge 2\rad(T) - 1$. Applying Theorem~\ref{thm:Diameter} to the tree $T$, we have that
\[
\gamma_k(G) = \gamma_k(T) \ge \frac{\diam(T)+1}{2k+1} \ge \frac{2\rad(T)}{2k+1} \ge \frac{2\rad(G)}{2k+1}. \hspace*{0.5cm} \Box
\]

\medskip
That the lower bound of Corollary~\ref{cor:Radius} is tight, may be seen by taking $G$ to be a path, $P_n$, of order~$n = 2\ell (2k+1)$ for some integer $\ell \ge 1$. Let $d = \diam(G)$ and let $r = \rad(G)$, and so $d = 2\ell (2k+1) - 1$ and $r = \ell (2k+1)$. In particular, we note that $d = 2r - 1$. By Proposition~\ref{prop:path}, $\gamma_k(G) = \frac{d + 1}{2k+1} = \frac{2r}{2k+1}$. As before by replacing each the internal vertices on the path with a clique of size at least~$\delta \ge 1$, we can obtain a graph with minimum degree~$\delta$ achieving the lower bound of Corollary~\ref{cor:Radius}.

We first prove a generalization of Theorem~\ref{thm:girth} by establishing a lower bound on the $k$-domination number of a graph in terms of its girth. We remark that when $k = 1$, Theorem~\ref{thm:Girth} is precisely Theorem~\ref{thm:girth}.

\begin{thm}
\label{thm:Girth}
For $k \ge 1$, if $G$ is a connected graph with girth $g$, then
\[
\gamma_{k}(G) \ge \frac{g}{2k+1}.
\]
\end{thm}
\proof The lower bound is trivial if $g \le 2k+1$. We may therefore assume that $g \ge 2k+2$, for otherwise the desired result is immediate. Let $C$ be a shortest cycle in $G$, and so $C$ has length~$g$. We note that the distance between two vertices in $V(C)$ is exactly the same in $C$ as in $G$.  We consider two cases, depending on the value of the girth.

\medskip \emph{Case 1. $2k+2 \le g \le 4k+2$.}  In this case, we need to show that $\gamma_{k}(G) \ge \lceil \frac{g}{2k+1} \rceil = 2$. Suppose, to the contrary, that $\gamma_{k}(G) = 1$. Then, $G$ contains a vertex $v$ that is within distance~$k$ from every vertex of $G$. In particular, $d(u,v) \le k$ for every vertex $u \in V(C)$. If $v \in V(C)$, then, since $C$ is a shortest cycle in $G$, we note that $d_C(u,v) = d_G(u,v) \le k$ for every vertex $u \in V(C)$. However, the lower bound condition on the girth, namely $g \ge 2k+2$, implies that no vertex on the cycle $C$ is within distance~$k$ in $C$ from every vertex of $C$, a contradiction. Therefore, $v \notin V(C)$.

By Lemma~\ref{Lem:CycleOutsider}, there exist two vertices $u,w \in V(C)$ such that a shortest $(v,u)$-path does not contain $w$ and a shortest $(v,w)$-path does not contain $u$. We show that we can choose $u$ and $w$ to be adjacent vertices on $C$. Let $w$ be a vertex of $C$ at maximum distance, say $d_w$, from~$v$ in $G$. Let $w_1$ and $w_2$ be the two neighbors of $w$ on the cycle $C$. If $d_G(v,w_1) = d_w$, then we can take $u = w_1$, and the desired property (that a shortest $(v,u)$-path does not contain $w$ and a shortest $(v,w)$-path does not contain $u$) holds. Hence, we may assume that $d_G(v,w_1) \ne d_w$. By our choice of the vertex $w$, we note that $d_G(v,w_1) \le d_w$, implying that $d_G(v,w_1) = d_w - 1$. Similarly, we may assume that $d_G(v,w_2) = d_w - 1$. Let $P_w$ be a shortest $(v,w)$-path. At most one of $w_1$ and $w_2$ belong to the path $P_w$. Renaming $w_1$ and $w_2$, if necessary, we may assume that $w_1$ does not belong to the path $P_w$. In this cas
 e, letting $u = w_1$ and letting $P_u$ be a shortest $(v,u)$-path, we note that $w \notin V(P_u)$. As observed earlier, $u \notin V(P_w)$. This shows that $u$ and $w$ can indeed be chosen to be neighbors on $C$.

Let $x$ be the last vertex in common with the $(v,u)$-path, $P_u$, and the $(v,w)$-path, $P_w$. Possibly, $x = v$. Then, the cycle obtained from the $(x,u)$-section of $P_u$ by proceeding along the edge $uw$ to $w$, and then following the $(w,x)$-section of $P_w$ back to~$x$, has length at most~$d_G(v,u) + 1 + d_G(v,w) \le 2k+1$, contradicting the fact that the girth $g \ge 2k+2$. Therefore, $\gamma_{k}(G) \ge 2$, as desired.

\medskip
\emph{Case 2. $g \ge 4k+3$.} Let $S$ be a minimum $k$-dominating set of $G$, and so $|S| = \gamma_k(G)$. Let $K = S \cap V(C)$ and let $L = S \setminus V(C)$. Thus, $S = K \cup L$. If $L = \emptyset$, then $S = K$ and the set $K$ is a $k$-dominating set of $C$, implying by Proposition~\ref{prop:cycle}, that $\gamma_k(G) = |S| = |K| \ge \gamma_k(C_g) = \lceil \frac{g}{2k+1} \rceil$, and the theorem holds. Hence we may assume that $|L| \ge 1$, for otherwise the desired result holds. We wish to show that $|K| + |L| = |S| \ge \lceil \frac{g}{2k+1} \rceil$. Suppose, to the contrary, that
\[
|K| \le \left\lceil \frac{g}{1+2k} \right\rceil -1 - |L|.
\]

As observed earlier, the distance between two vertices in $V(C)$ is exactly the same in $C$ as in $G$. This implies that each vertex of $K$ (recall that $K \subseteq V(C)$) is within distance~$k$ from exactly~$2k+1$ vertices of $C$. Thus, the set $K$ $k$-dominates at most
\[
\begin{array}{lcl}
|K|(2k+1) & \le & \left( \left\lceil \frac{g}{2k+1} \right\rceil -1 - |L| \right) (2k+1) \1 \\
& \le & \left( \frac{g+2k}{2k+1}  - 1 - |L| \right) (2k+1) \1 \\
& = & g -  1 - |L|(2k+1) \1
\end{array}
\]
vertices from $C$. Consequently, since $|C(V)| = g$, there are at least $|L|(2k+1) + 1$ vertices of $C$ which are not $k$-dominated by vertices of $K$, and therefore must be $k$-dominated by vertices from $L$.  Thus, by the Pigeonhole Principle, there is at least one vertex, call it $v$, in $L$ that $k$-dominates at least $2k+2$ vertices in $C$. By Lemma~\ref{Lem:CycleOutsider},  there exist two vertices $u,w \in V(C)$ that are both $k$-dominated by $v$ and such that a shortest $(u,v)$-path, $P_u$ say, (from $u$ to $v$) does not contain $w$ and a shortest $(w,v)$-path, $P_w$ say, (from $w$ to $v$) does not contain $u$. Analogously as in the proof of Lemma~\ref{Lem:CycleOutsider}, we can choose the vertex $u$ to be a vertex of $C$ at minimum distance from~$v$ in $G$. Thus, the vertex $u$ is the only vertex on the cycle $C$ that belongs to the path~$P_u$. Combining the paths $P_u$ and $P_w$ produces a $(u,w)$-walk of length at most~$d_G(u,v) + d_G(v,w) \le 2k$, implying that $d
 _G(u,w) \le 2k$. Since $C$ is a shortest cycle in $G$, we therefore have that $d_C(u,w) = d_G(u,w) \le 2k$. The cycle $C$ yields two $(w,u)$-paths. Let $P_{wu}$ be the $(w,u)$-path on the cycle $C$ of shorter length  (starting at $w$ and ending at $u$). Thus, $P_{wu}$ has length~$d_C(u,w) \le 2k$. Note that the path $P_{wu}$ belongs entirely on the cycle $C$. Let $x \in V(C)$ be the last vertex in common with the $(w,v)$-path, $P_w$, and the $(w,u)$-path, $P_{wu}$. Possibly, $x = w$. However, note that $x \ne u$ since $u \notin V(P_w)$. Let $y$ be the first vertex in common with the $(x,v)$-subsection of the path $P_w$ and with the $(u,v)$-path $P_u$. Possibly, $y = v$. However, note that $y \ne x$ since $x \notin V(P_u)$ and $V(P_u) \cap V(C) = \{u\}$. Using the $(x,u)$-subsection of the path $P_{wu}$, the $(x,y)$-subsection of the path $P_w$, and the $(u,y)$-subsection of the path $P_u$ produces a cycle in $G$ of length at most~$d_G(u,v) + d_G(w,v) + d_G(u,w) \le k + k + 2
 k = 4k$, contradicting the fact that the girth $g \ge 4k + 3$. Therefore, $\gamma_k(G) = |S| = |K| + |L| \ge \lceil \frac{g}{2k+1} \rceil$, as desired.~\qed


\newpage
\section{Direct Product Graphs}

The \emph{direct product graph}, $G \times H$, of graphs $G$ and $H$ is the graph with vertex set $V(G) \times V(H)$ and with edges $(g_1,h_1)(g_2,h_2)$, where $g_1g_2\in E(G)$ and $h_1h_2\in E(H)$. Let $A \subseteq V(G \times H)$. The \emph{projection of $A$} onto $G$ is defined as
\begin{equation*}
P_G(A) = \{g\in V(G) \colon (g,h)\in A \:\: \text{for some} \:\: h\in V(H)\}.
\end{equation*}
Similarly, the projection of $A$ onto $H$ is defined as
\begin{equation*}
P_H(A)=\{g\in V(H) \colon (g,h)\in A\:\:\text{for some}\:\:h\in V(G)\}.
\end{equation*}

For a detailed discussion on direct product graphs, we refer the reader to the handbook on graph products~\cite{ProductGraphs}. There have been various studies on the domination number of direct product graphs. For example, Meki\v{s}~\cite{Mekis} proved the following lower bound on the domination number of direct product graphs. Recall that for every graph $G$, $\gamma(G) = \gamma_1(G)$.

\begin{thm}{\rm (\cite{Mekis})}
If $G$ and $H$ are connected graphs, then
\[
\gamma(G\times H)\ge \gamma(G) + \gamma(H) -1.
\]
\label{Mekis}
\end{thm}

\vskip -0.5cm
Staying within the theme of our previous results, we now prove a projection lemma which will enable us generalize the result of Theorem~\ref{Mekis} on the domination number to the $k$-domination number.

\begin{lem}{\rm (Projection Lemma)}
Let $G$ and $H$ be connected graphs. If $D$ is a $k$-dominating set of $G \times H$, then $P_G(D)$ is a $k$-dominating set of $G$ and $P_H(D)$ is a $k$-dominating set of $H$.
\label{Projection Lemma}
\end{lem}
\proof Let $D \subseteq V(G\times H)$ be a $k$-dominating set of $G \times H$. We show firstly that $P_G(D)$ is a $k$-dominating set of $G$. Let $g$ be a vertex in $V(G)$.  If $g \in P_G(D)$, then $g$
is clearly $k$-dominated by $P_G(D)$. Hence, we may assume that $g \in V(G) \setminus P_G(D)$. Let $h$ be an arbitrary vertex in $V(H)$. Since $g \notin P_G(D)$, the vertex $(g,h) \notin D$. However, the set $D$ is a $k$-dominating set of $G \times H$, and so $(g,h)$ is within distance~$k$ from $D$ in $G$; that is, $d_{G \times H}((g,h),D) \le k$. Let $(g_0,h_0), (g_1,h_1), \ldots, (g_r,h_r)$ be a shortest path from $(g,h)$ to $D$ in $G \times H$, where $(g,h) = (g_0,h_0)$ and $(g_r,h_r) \in D$. By assumption, $1 \le r \le k$. For $i \in \{0,\ldots,r-1\}$, the vertices $(g_i,h_i)$ and $(g_{i+1},h_{i+1})$ are adjacent in $G \times H$. Hence, by the definition of the direct product graph, the vertices $g_i$ and $g_{i+1}$ are adjacent in $G$, implying that $g_0g_1 \ldots g_r$ is a $(g_0,g_r)$-walk in $G$ of length~$r$. This in turn implies that there is a $(g_0,g_r)$-path in $G$ of length~$r$. Recall that $g = g_0$ and $1 \le r \le k$. Since $(g_r,h_r) \in D$, the vertex $g_r \i
 n P_G(D)$. Hence, there is a path from $g$ to a vertex of $P_G(D)$ in $G$ of length at most~$k$. Since $g$ is an arbitrary vertex in $V(G)$, the set $P_G(D)$ is therefore a $k$-dominating set of $G$. Analogously, the set $P_H(D)$ is a $k$-dominating set of $H$.~\qed

\medskip
Using our Projection Lemma, we are now in a position to generalize Theorem~\ref{Mekis}.

\begin{thm}
If $G$ and $H$ are connected graphs, then
\[
\gamma_k(G\times H) \ge \gamma_k(G) + \gamma_k(H) -1.
\]
\label{thm:Direct}
\end{thm}
\vskip -0.5cm
\noindent
\proof
Let $D \subseteq V(G \times H)$ be a minimum $k$-dominating set of $G \times H$. Suppose, to the contrary, that $|D| \le \gamma_k(G) + \gamma_k(H) - 2$. We will refer to this supposition as $(*)$.
By Lemma~\ref{Projection Lemma}, $P_G(D)$ is a $k$-dominating set of $G$ and $P_H(D)$ is a $k$-dominating set of~$H$. Therefore, we have that $|D| \ge |P_G(D)| \ge \gamma_k(G)$ and $|D| \ge |P_H(D)| \ge \gamma_k(H)$.
If $\gamma_k(G) = 1$, then, by $(*)$,
$\gamma_k(H) - 1 \ge |D| \ge \gamma_k(H)$,
a contradiction.  Therefore, $\gamma_k(G) \ge 2$. Analogously, $\gamma_k(H) \ge 2$. Recall that $|P_G(D)| \ge \gamma_k(G)$. We now remove vertices from the set $P_G(D)$ until we obtain a set, $D_G$ say, of cardinality exactly~$\gamma_k(G) - 1$. Thus, $D_G$ is a proper subset of $P_G(D)$ of cardinality~$\gamma_k(G) - 1$. Since $D_G$ is not a $k$-dominating set of $G$, there exists a vertex $g \in V(G)$ that is not $k$-dominated by the set $D_G$ in $G$; that is, $d_G(g,D_G) > k$.
Let $D_G = \{g_1,\ldots,g_t\}$, where $t = \gamma_k(G) - 1 \ge 1$. For each $i \in [t]$, there exists a (not necessarily unique) vertex $h_i \in V(H)$ such that $(g_i,h_i) \in D$ (since $D_G = P_G(D)$). We now consider the set
\[
D_0 = \{(g_1,h_1),\dots, (g_t,h_t)\},
\]
and note that $D_0 \subset D$ and $|D_0| = \gamma_k(G) - 1$. By $(*)$, we note that

\[
\begin{array}{lcl}
|P_H(D \setminus D_0)| & \le & |D \setminus D_0| \1 \\
& = & |D| - |D_0| \1 \\
& \le & (\gamma_k(G) + \gamma_k(H) - 2) - (\gamma_k(G) - 1) \1 \\
& = &  \gamma_k(H) -1 \1 \\
& < &  \gamma_k(H).
\end{array}
\]

Hence, there exists a vertex $h \in V(H)$ that is not $k$-dominated by the set $P_H(D \setminus D_0)$ in $H$; that is, $d_H(h,P_H(D \setminus D_0)) > k$. We now consider the vertex $(g,h) \in V(G \times H)$. Since $D$ is a $k$-dominating set of $G \times H$, the vertex $(g,h)$ is $k$-dominated by some vertex, say $(g^*,h^*)$, of $D$ in $G \times H$. An analogous proof as in the proof of Lemma~\ref{Projection Lemma} shows that $d_G(g,g^*) \le k$ and $d_H(h,h^*) \le k$. If $(g^*,h^*) \in D \setminus D_0$, then $h^* \in P_H(D \setminus D_0)$, implying that $d_H(h,P_H(D \setminus D_0)) \le d_H(h,h^*) \le k$, a contradiction. Hence, $(g^*,h^*) \in D_0$. This in turn implies that $g^* \in P_G(D_0) = G_D$. Thus, $d_G(g,D_G) \le d_G(g,g^*) \le k$, contradicting the fact that $d_G(g,D_G) > k$. Therefore, the supposition that $|D| = \gamma_k(G) + \gamma_k(H) - 2$ must be false, and
the result follows.~\qed


\end{document}